\documentclass{article}
\usepackage{amsmath,amssymb,amsfonts,amsthm}
\usepackage{graphicx}
 \usepackage{hyperref}
\makeatletter
\renewcommand\section{\@startsection {section}{1}{\z@}%
                                   {-3.5ex \@plus -1ex \@minus -.2ex}%
                                   {2.3ex \@plus.2ex}%
                                   {\normalfont\Large\bfseries
                                    \setcounter{equation}{0}
                                    \setcounter{theorem}{0}}}
\makeatother

\newtheorem{theorem}{Theorem}

\theoremstyle{definition}
\newtheorem{thm}[theorem]{Theorem}
\newtheorem{rem}[theorem]{Remark}

\newtheorem{defn}[theorem]{Definition}

\newtheorem{lem}[theorem]{Lemma}
\newlength{\defbaselineskip}
\newcommand{\setlinespacing}[1]%
           {\setlength{\baselineskip}{#1 \defbaselineskip}}
          
 \newcommand{\qbinom}[2] {{#1 \brack #2}_q}
  \newcommand{\pbinom}[2] {{#1 \brack #2}_p}

\newcommand{\set}[1]{\left\{#1\right\}}

\newcommand{\be }{\begin{equation}}
\newcommand{\ee }{\end{equation}}
\newcommand{\bbb }{\begin{eqnarray}}
\newcommand{\eee }{\end{eqnarray}}
\newcommand{\bea }{\begin{eqnarray*}}
\newcommand{\eea }{\end{eqnarray*}}
\newcommand{\ed }{\end{document}}

\theoremstyle{plain}

\usepackage{amsfonts}

\newtheorem{exa}{Example}
\newtheorem{prop}{Proposition}

\date{ }
\title{ $q$-type Lidstone expansions and an  interpolation problem\\ for entire functions}

\author{\small Mourad E.H.I. Ismail  $^{\lowercase{a}}$, Zeinab S. I. Mansour $^{\lowercase{b}}$ \\
\small$^{\lowercase{a}}$ Department of Mathematics,
University of Central Florida, Orlando, Florida, USA. \\
\small E-mail: mourad.eh.ismail@gmail.com\\
\small $^{\lowercase{b}}$
Department of Mathematics,
\small Faculty of Science, Cairo University, Giza, Egypt.\\
\small E-mail:  zeinab@sci.cu.edu.eg\\}

\begin{document}
\thispagestyle{empty} \maketitle

\begin{abstract}
In this paper, we expand functions of specific  $q$-exponential growth in terms  of  its even  (odd) Askey- Wilson $q$-derivatives  at $0$ and $\eta=(q^{1/4}+q^{-1/4})/2$. This expansion is a $q$-version of the celebrated Lidstone expansion theorem, where 
we expand the function in $q$-analogs of Lidstone polynomials, i.e.,   q-Bernoulli
 and $q$-Euler polynomials as in the classical case.  We also raise and solve a $q$-extension of the problem of representing an entire function of the form 
 $f(z)=g(z+1)-g(z)$, where $g(z)$ is also an entire function of the same order as $f(z)$. 
\end{abstract}

 \noindent {\it 2020 Mathematics Subject Classification:} 05A30, 41A58, 39A13, 30E20 \\
\noindent {\it Keywords:}{  $q$- Lidstone expansion theorem,  Askey-Wilson divided difference operator, 
$q$-Bernoulli and Euler  polynomials. }

\section{Introduction and Preliminaries }

The Taylor series expands  an analytic function $f$ 
with prescribed values of $f$ and all its derivatives at a single 
point. 
The Lidstone expansion refers to the expansion of an entire 
function in power series at  $z-a$ and $z-b$ when
 $f^{(2n)}(a)$ and $f^{(2n)}(b)$ are prescribed. 

One question is what happens if,  instead of the derivative 
operator,  we use other operators. For example,  the $q$-Taylor series using 
the $q$-difference operator was developed  by F. H. Jackson~\cite{Jackson} 
and more rigorously by Annaby and Mansour 
\cite{Annaby-Mansour}. The $q$-Taylor series expansion 
in terms of the Askey--Wilson operator was developed by 
Ismail and Stanton \cite{Ismail-Stanton-2003,Ismail-Stanton-2003-2}. These are 
instances of the general interpolation problem using divided 
difference operators,  see,  e.g. \cite{Guelfond}.

For convenience,  we now define the operators of interest in this work.  The action of the $q$-difference operator $D_q$  on a function $f$ is $D_qf(x)=\frac{f(x)-f(qx)}{x-qx}$, $x\neq 0$ and the $q$-derivative at zero is usually defined as $f'(0)$ if $f$ is differentiable at zero. The Askey-Wilson divided difference operator $\mathcal{D}_q$, \cite{Askey-Wilson-1985}, is  defined as follows. 
Given a polynomial $f$, we set $\breve{f}(e^{i\theta}):=f(x)$,  $x=\cos\theta$ . In this notation, the Askey-Wilson divided difference operator  $\mathcal{D}_q$ is defined by 
\be 
\mathcal{D}_q f(x)=\dfrac{\breve{f}(q^{1/2}z)-\breve{f}(q^{-1/2}/z)}{i(q^{1/2}-q^{-1/2})\sin\theta},\;x= (z+1/z)/2. 
\ee

We now state the Lidstone expansion  theorem.  It asserts  that an entire function $f(z)$ of exponential type less than $\pi$ has a convergent (Lidstone) representation as
\[f(z)=\sum_{n=0}^{\infty}A_n(1-z) f^{(2n)}(0)+\sum_{n=0}^{\infty}A_n(z) f^{(2n)}(1),\]
where the polynomial $A_n(z)$ is called Lidstone polynomials. 
See~\cite{Boas-Buck}. It is known that $A_n(z)$ is a constant 
multiplier of the Bernoulli polynomial of order $2n+1$. More 
Precisely, 
\[A_n(z)=\frac{2^{2n+1}}{(2n+1)!} B_{2n+1}(\frac{z+1}{2}).\]
See~\cite[Eq. (3.5)]{Whittaker}.  Ismail and Mansour in 
\cite{Ismail-Mansour} introduced a $q$ analog of Lidstone 
theorem where they expand a  class of entire functions of 
$q$-exponential growth in terms  of Jackson $q$ derivatives of 
even degree at 0 and 1.  See also~\cite{Mansour-Towailb1, 
Mansour-Towailb2} for applications to the $q$-Lidstone theorem  we 
 introduced.
 This paper aims to introduce  $q$-Lidstone expansion theorems 
 for the Askey-Wilson $q$-difference operators. We also 
 introduce a $q$-version of the problem raised by 
 Guichard~\cite{Guichard}  about the existence of  an entire 
 function $g(z)$ such that an entire function $f(z)$ can be 
 represented  of the form $g(z+1)-g(z)=f(z)$.

  The rest of this section contains definitions and  
 preliminary results used  in the 
 sequel. In Section 2, we introduce a new $q$-analog of 
 Bernoulli and Euler polynomials and study their properties. 
 Section 3 contains the $q$-Lidstone expansion theorem 
 associated with the Askey-Wilson divided difference operator.   Finally, 
 In section 4, we formulate  and solve  a $q$-analog of the 
 interpolation problem introduced by  Guichard in
  \cite{Guichard}
 and developed by Whittaker~ \cite{Whittaker}.

Throughout this paper, unless otherwise is stated, $q$ is a positive number less than 1, $\mathbb{N}$ is the set of positive integers, and $\mathbb{N}_0$ is the set of nonnegative integers. We also follow Gasper and Rahman~\cite{GR} for the definitions of $q$-shifted factorials, $q$-binomial coefficients, $q$-numbers, and Jackson $q$-Bessel functions. 
Ismail and Zhang, see~ \cite{Ismail-Zhang-94},  introduced the  $q$-exponential function 
\be\begin{gathered}
\mathcal{E}_q(\cos\theta, \cos\phi;\alpha):=\\\frac{(\alpha^2;q^2)_{\infty}}{(q\alpha^2;q^2)_{\infty}}\sum_{n=0}^{\infty}\frac{(\alpha e^{-i\phi})^n}{(q;q)_n} q^{n^2/4}\left(-e^{i(\phi+\theta)} q^{(1-n)/2},-e^{i(\phi-\theta)} q^{(1-n)/2};q\right)_n.
\end{gathered} 
 \ee

The function $\mathcal{E}_q(x;\alpha)$ is defined as
\be 
\mathcal{E}_q(x;\alpha):=\mathcal{E}_q(\cos\theta, 0;\alpha),\; x=\cos \theta.
\ee
In other words,
\be
\mathcal{E}_q(\cos\theta;\alpha):=\frac{(\alpha^2;q^2)_{\infty}}{(q\alpha^2;q^2)_{\infty}}\sum_{n=0}^{\infty}
\frac{(-i\alpha)^n}{(q;q)_n} q^{n^2/4}\left(-ie^{i\theta} q^{(1-n)/2},-ie^{-i\theta} q^{(1-n)/2};q\right)_n.
 \ee
The $q$-exponential function $\mathcal{E}_q(x;\alpha)$  has the properties \[\mathcal{E}_q(0;\alpha)=1\;\mbox{ and}\; 
 \lim_{q\to 1} \mathcal{E}_q(x; \frac{(1-q)t}{2})=\exp(xt).\]    See~\cite[Lemma 14.1.1]{Ismail-book}. In \cite{Berg-Pedersen}, 
 
Ramis~\cite{Ramis} defined an entire function $f$ to have a $q$-exponential growth of order $k$ and a finite type if there exist real numbers $K$, $\alpha$, $K>0$, such that 
\be\label{Ramis} |f(x)|< K |x|^{\alpha}\exp\left(\frac{k\ln^2|x|}{2\ln^2 q}\right).\ee
In general, for an entire function, $f$, of order zero,  the logarithmic order, $\rho$,  is  defined by 
\[\rho=\inf\{a>0:\log M(r;f)\leq   (\log r)^a , r\geq r_0\}.\] 
If $\rho <\infty$, the logarithmic type is 
\[\tau=\inf\{\beta >0:\log M(r;f)\leq \beta (\log r)^a, r\geq r_0\}.\]
See \cite{Berg-Pedersen}.  Therefore,  if $f$ is a function of $q$-exponential growth of order $k$, then   $f$ has a logarithmic order two and type $\frac{k}{2\ln^2 q}$. 
From \cite[Lemma 2.2]{Ramis}, if  $f$ satisfies the condition \eqref{Ramis} and $f(x)=\sum_{n=0}^{\infty}a_n x^n$ then
\be\label{Coeff.} |a_n|\leq K e^{-(n-\alpha)^2\frac{\ln^2 q}{2k}}\quad (n\in\mathbb{N}).\ee
Ismail and Stanton proved the following theorem, see \cite{Ismail-Stanton-2003}. 
\begin{thm}
For $|y|<1$, the maximum modulus of $\mathcal{E}_q(\cdot;y)$ has the property that
\be
\lim_{r\to\infty} \sup \frac{ln M(r;\mathcal{E}_q)}{\ln^2 r}=\frac{1}{\ln q^{-1}},
\ee
where \[M(r;f):=\sup\{|f(z)|: |z|=r\}.\]
\end{thm}
It is worth noting that $\mathcal{E}_q(x;y)$ is  an entire function in $x$ only for $|y|<1$ but $(qy^2;q^2)_{\infty}\mathcal{E}_q(x;y)$ is an  entire function in $x$ for all $y\in\mathbb{C}$. See \cite{Ismail-book, Suslov-2003}. 
Also, $\mathcal{E}_q(x;t)$   appears in the generating function of the continuous $q$-Hermite polynomials  as 
\be\label{Hermite-rep}
(q^2 t^2;q^2)_{\infty}\mathcal{E}_q(x;t)=\sum_{n=0}^{\infty} \frac{q^{n^2/4}}{(q;q)_n}t^n H_n(x|q).
\ee
Moreover, 
\be\label{rho-rep}
\mathcal{E}_q(x;t)=\sum_{n=0}^{\infty} \frac{q^{n^2/4}}{(q;q)_n}t^n \rho_n(x),
\ee
where 
\[\rho_n(x):=(1+e^{2i\theta}) (-q^{2-n}e^{2i\theta};q^2)_{n-1} e^{-in\theta}\quad (n\geq 1)\]
and $\rho_0(x)=1$. 
One can verify that 
\[\mathcal{D}_q\rho_n(x)= 2q^{\frac{1-n}{2}}\frac{1-q^n}{1-q}\rho_{n-1}(x).\]
See~\cite[Theorem 14.2.1, (14.6.7))]{Ismail-book}. 
Recall that the continuous $q$-Hermite polynomials defined by the generating function
\be
\sum_{n=0}^{\infty}H_n(\cos\theta|q)\frac{t^n}{(q;q)_n}=\frac{1}{(te^{i\theta}, te^{-i\theta};q)_{\infty}},
\ee
which yields 
\be  
H_n(cos\theta|q)=\sum_{k=0}^{n}\frac{(q;q)_n}{(q;q)_k (q;q)_{n-k}} \cos(n-2k)\theta.
\ee

 The basic cosine and sine functions are denoted by  $C_q(x,y;w)$ and   $S_q(x,y;w)$, respectively,   and defined  through the identity 
 \[\mathcal{E}_q(x,y;iw)=C_q(x,y;w)+iS_q(x,y;w).\]
 They are $q$-analogs of $\cos w(x+y)$ and $\sin w(x+y)$, respectively. In the sense that 
 
 \[\begin{split}
 \lim_{q\to 1^-}C_q(x,y;\frac{w(1-q)}{2})=\cos w(x+y),
 \lim_{q\to 1^-} S_q(x,y;\frac{w(1-q)}{2})=\sin w(x+y).
 \end{split}
 \]
 The special case $y=0$ ($\phi=\frac{\pi}{2}$) are
 
 \be\label{Eq:BC}\begin{split}
 C_q(x;w)&=\dfrac{(-w^2;q^2)_{\infty}}{(-qw^2;q^2)_{\infty}}{}_2\phi_1\left(-qe^{2i\theta},-qe^{-2i\theta};q;q^2;-w^2\right)\\
 &=\dfrac{(-w^2;q^2)_{\infty}}{(-qw^2;q^2)_{\infty}}\sum_{k=0}^{\infty}\dfrac{(-qe^{2i\theta},-qe^{-2i\theta};q^2)_k}{(q;q)_{2k}}(-w^2)^k
 \end{split}\ee
 
\be\label{Eq:BS} \begin{split}
 S_q(x;w)&=\dfrac{(-w^2;q^2)_{\infty}}{(-qw^2;q^2)_{\infty}}\frac{2q^{1/4}w}{1-q}x\;{}_2\phi_1\left(-q^2e^{2i\theta},-q^2e^{-2i\theta};q^3;q^2;-w^2\right)\\
 &=2q^{1/4} x\dfrac{(-w^2;q^2)_{\infty}}{(-qw^2;q^2)_{\infty}}\sum_{k=0}^{\infty}\dfrac{(-q^2e^{2i\theta},-q^2e^{-2i\theta};q^2)_k}{(q;q)_{2k+1}}(-1)^k w^{2k+1}.
 \end{split}\ee
 Since 
 \[\mathcal{D}_q\mathcal{E}_q(x;w)=\frac{2q^{1/4}}{(1-q) } w\mathcal{E}_q(x;w),\]
then 

 \[\mathcal{D}_q S_q(x;w)=\frac{2q^{1/4}}{(1-q) } wC_q(x;w),\;\mathcal{D}_q C_q(x;w)=-\frac{2q^{1/4}}{(1-q) } S_q(x;w).\]
 
\begin{prop} The basic sine and cosine functions have the series representation:
\bea
 S_q(x;w)&=&q^{1/4}\sum_{n=0}^{\infty}(-1)^n q^{n^2+n}\frac{w^{2n+1}}{(q;q)_{2n+1}} \rho_{2n+1}(x),\\
 C_q(x;w)&=&\sum_{n=0}^{\infty}(-1)^n q^{n^2}\frac{w^{2n}}{(q;q)_{2n}} \rho_{2n}(x).
 \eea
 \end{prop}
 \begin{proof}
 The proof follows from the identity 
 \[\mathcal{E}_q(x;iw)=C_q(x;w)+i S_q(x;w)\]
 and the identity \eqref{rho-rep}.
 \end{proof}
  The basic sine and cosine functions are entire functions in $x$ when $|w|<1$. Ismail and Stanton
   proved that when 
  $\eta:=(q^{1/4}+q^{-1/4})/2$ 
  \[\begin{split}
  \mathcal{E}_q(\eta;w)&=\dfrac{(-w;q^{1/2})_{\infty}}{(qw^2;q^2)_{\infty}}=\dfrac{(-w;q)_{\infty}}{(q^{1/2}w;q)_{\infty}},\\
  S_q(\eta;w)&=\dfrac{(-iw;q^{1/2})_{\infty}-(iw;q^{1/2})_{\infty}}{2i(-qw^2;q^2)_{\infty}}=\sum_{k=0}^{\infty}(-1)^k\dfrac{(-q^{-1/2};q)_{2k+1}}{(q;q)_{2k+1}}(q^{1/2}w)^{2k+1},\\
  C_q(\eta;w)&=\dfrac{(-iw;q^{1/2})_{\infty}+(iw;q^{1/2})_{\infty}}{2(-qw^2;q^2)_{\infty}}=\sum_{k=0}^{\infty}(-1)^k\dfrac{(-q^{-1/2};q)_{2k}}{(q;q)_{2k}}(q^{1/2}w)^{2k}.
  \end{split}
  \]
 Moreover, 
 \be \label{Eq:5}\begin{split}
 S_q(\eta;w)&=\dfrac{(q;q)_{\infty}}{(q^{1/2};q)_{\infty}} \frac{w^{1/2}}{(-qw^2;q^2)_{\infty}} J_{1/2}^{(2)}(2w;q),\\
 C_q(\eta;w)&=\dfrac{(q;q)_{\infty}}{(q^{1/2};q)_{\infty}} \frac{w^{1/2}}{(-qw^2;q^2)_{\infty}} J_{-1/2}^{(2)}(2w;q).
\end{split} 
 \ee
In~\cite{Ism1},  Ismail
proved that the zeros of $J_{\nu}^{(2)}(z;q)$ are real and simple
and  the zeros of $z^{-\nu}J_{\nu}^{(2)}(z;q)$ and
$z^{-\nu}J_{\nu+1}^{(2)}(z;q)$ interlace.  Ismail, cf.~\cite{Ism1},
also conjectured that if $\set{z_{m,\nu}}_{m=1}^{\infty}$ are the
positive zeros of $J_{\nu}^{(2)}(z;q)$, then
\begin{equation}\label{Ism}
z_{m,\nu}=c_1q^{-m/\alpha}+c_2q^{\beta
m}\left(1+O(q^{\gamma\/m})\right),\quad\mbox{as}\quad m\to\infty,
\end{equation}
for some constants $c_1,\,c_2,\,\alpha,\,\beta$, and  $\gamma$.
Hayman proved in~\cite{Hay4}
that, for arbitrary $k$,  the positive zeros
$\set{z_{m,\nu}}_{m=1}^{\infty}$ of $J_{\nu}^{(2)}(z;q)$ have the
following asymptotic expansion
\begin{equation}\label{hay1}
z_{m,\nu}=2q^{-m}q^{\frac{-\nu+1}{2}}\left\{1+\sum_{n=1}^{k}b_{n}q^{mn}+O(q^{(k+1)m})\right\},
\end{equation}
for sufficiently large $m$, where the constants $b_n$,
$n=1,2,\ldots,k$ depend on $q$ and $\nu$ and can be computed
iteratively. In particular,  when $k=1$ we have
\begin{equation}\label{hay2}
z_{m,\nu}=2q^{-m}q^{\frac{-\nu+1}{2}}\big(1+O(q^{m})\big)\quad (m\to\infty).
\end{equation}
 Therefore the zeros
$\set{\xi_{m,\nu}}_{m=1}^{\infty}$ of $J_{\nu}^{(2)}(z;q^2)$ have
the asymptotic behavior
\begin{equation}\label{hay3}
\xi_{m,\nu}=2q^{-2m}q^{-\nu+1}\big(1+O(q^{2m})\big),
\end{equation}
for sufficiently  large $m$.

 \begin{prop}\label{prop:1}
 The functions $S_q(\eta,w)$ and $C_q(\eta,w)$ are analytic for $|w|<q^{-1/2}$. Moreover, if $w_1 $ is the smallest positive zero of $S_q(\eta,w)$, then \[w_1^2\geq q^{-3/2}(1-q)(1-q^{3/2}).\]
 \end{prop}
 \begin{proof}
 The proof of the first part of the proposition follows at once  from \eqref{Eq:5}. So, we only prove  the second part of the proposition. From \eqref{Eq:5},  
 \[(-qw^2;q^2)_{\infty}S_q(\eta;w)=\sum_{k=0}^{\infty}(-1)^k\frac{q^{k(k+\frac{1}{2})}}{(q^{1/2};q^{1/2})_{2k+1}} w^{2k+1}. \]
 Set $a_k:=\dfrac{q^{k(k+\frac{1}{2})}}{(q^{1/2};q^{1/2})_{2k+1}} w^{2k+1} $, then \[\frac{a_{k+1}}{a_k}=\frac{q^{2k+\frac{3}{2}}w^2}{(1-q^{k+1})(1-q^{k+\frac{3}{2}})}<\frac{q^{\frac{3}{2}}w^2}{(1-q)(1-q^{\frac{3}{2}})}<1 \]
 for all $k\in \mathbb{N}_0$ if $|w|^2<q^{-3/2}(1-q)(1-q^{3/2})$.  Consequently, 
 \[(-qw^2;q^2)_{\infty}S_q(\eta,w)=\sum_{k=0}^{\infty}\left(a_{2k}-a_{2k+1}\right)>0\]  if $|w|^2<q^{-3/2}(1-q)(1-q^{3/2})$. Therefore,
  $w_1^2\geq q^{-3/2}(1-q)(1-q^{3/2})$.
 \end{proof} 
Proposition 1 shows that $w_1$ may be less than 1 or greater than 1 based  on the value of $q\in(0,1)$. 

\section{ $q$-analogs of the Bernoulli and Euler numbers}
Recall that   the generating functions define the classical Bernoulli and Euler's polynomials are 
\be\label{B-E}
\frac{w e^{xw}}{e^w-1}=\sum_{n=0}^{\infty}B_n(x)\frac{w^n}{n!},\;\mbox{and}\;
\frac{2 e^{xw}}{e^w+1}=\sum_{n=0}^{\infty}E_n(x)\frac{w^n}{n!},
\ee
respectively.   The Bernoulli numbers $(\beta_n)_n$ are defined as the values of the Bernoulli polynomials at $x=0$. I.e.
\be
\frac{w }{e^w-1}=\sum_{n=0}^{\infty}\beta_n\frac{w^n}{n!}.
\ee
Suslov, see ~\cite{Suslov-2003}, introduced    $q$-analogs of the Bernoulli  and Euler polynomials through the generating functions

\be \label{S:def}
\begin{split}
\sum_{n=0}^{\infty}\mathcal{B}_n(x;q) w^n&=\dfrac{w \mathcal{E}_q(x;w)}{\mathcal{E}_q(\eta;w)-\mathcal{E}_q(\eta;-w)}\\
&=\dfrac{w (qw^2;q^2)_{\infty}\mathcal{E}_q(x;w)}{(-w;q^{1/2})_{\infty}-(w;q^{1/2})_{\infty}},
\end{split}\ee

\be \label{S:ENdef}
\begin{split}
\sum_{n=0}^{\infty}\mathcal{E}_n(x;q) w^n&=\dfrac{ \mathcal{E}_q(x;w)}{\mathcal{E}_q(\eta;w)+\mathcal{E}_q(\eta;-w)}\\
&=\dfrac{ (qw^2;q^2)_{\infty}\mathcal{E}_q(x;w)}{(-w;q^{1/2})_{\infty}+(w;q^{1/2})_{\infty}},
\end{split}\ee
respectively. Suslov also defined  the basic Bernoulli numbers $\mathcal{B}_n(q)$ as the value of the  basic Bernoulli polynomials at $-\eta$,  and  the Euler's number as the value of Euler's polynomial at zero.
I.e. 
\be
 \label{S:def2}
\dfrac{w (w;q^{1/2})_{\infty}}{(-w;q^{1/2})_{\infty}-(w;q^{1/2})_{\infty}}=\sum_{n=0}^{\infty}\mathcal{B}_n(q) w^n ,
\ee
\be \label{S:ENdef2}
\dfrac{ (w;q^{1/2})_{\infty}}{(-w;q^{1/2})_{\infty}+(w;q^{1/2})_{\infty}}=\sum_{n=0}^{\infty}\mathcal{E}_n(q) w^n ,
\ee
\[ \mathcal{B}_n(q)=\frac{2^{n-1}}{(q;q^{1/2})_{n-1}}B_n(q^{1/2}).\]

\begin{defn}

We define new $q$-analogs of Bernoulli  and Euler's polynomial by
\be\label{IM:def} \begin{split}
\sum_{n=0}^{\infty}\beta_n(x;q) w^n&=\dfrac{w\mathcal{E}_q(\eta;-w) \mathcal{E}_q(x;w)}{\mathcal{E}_q(\eta;w)-\mathcal{E}_q(\eta;-w)}\\
&=\dfrac{w (w;q^{1/2})_{\infty}\mathcal{E}_q(x;w)}{(-w;q^{1/2})_{\infty}-(w;q^{1/2})_{\infty}},
\end{split}
\ee

\be \begin{split}\label{IM:Edef}
\sum_{n=0}^{\infty}{E}_n(x;q) w^n&=\dfrac{2\mathcal{E}_q(\eta;-w) \mathcal{E}_q(x;w)}{\mathcal{E}_q(\eta;w)+\mathcal{E}_q(\eta;-w)}\\
&=
\dfrac{2 (w;q^{1/2})_{\infty}\mathcal{E}_q(x;w)} {(-w;q^{1/2})_{\infty}+(w;q^{1/2})_{\infty}}.
\end{split}
\ee
\end{defn}

We define the $q$-Bernoulli number $(\beta_n(q))_n$ as  the $q$-Bernoulli polynomials' values at $x=0$. I.e.

\be \label{BND2}
\dfrac{w\mathcal{E}_q(\eta; -w)}{\mathcal{E}_q(\eta; w)-\mathcal{E}_q(\eta; -w)}=\sum_{n=0}^{\infty}\beta_n(q) w^n.
\ee

\begin{prop} For $n\in\mathbb{N}$
\[ \mathcal{D}_q \mathcal{B}_n(x;q)=\frac{2 q^{1/4}}{1-q}\mathcal{B}_{n-1}(x;q),\;  \mathcal{D}_q \beta_n(x;q)=\frac{2 q^{1/4}}{1-q}\beta_{n-1}(x;q),\]
\[      \mathcal{D}_q\mathcal {E}_n(x;q)=\frac{2 q^{1/4}}{1-q}\mathcal{E}_{n-1}(x;q),\;  \mathcal{D}_q \widetilde{E}_n(x;q)=\frac{2 q^{1/4}}{1-q}\widetilde{E}_{n-1}(x;q).\]
\end{prop}
\begin{proof} Let $n$ be a positive integer. 
We prove 
$ \mathcal{D}_q \beta_n(x;q)=\frac{2 q^{1/4}}{1-q}\beta_{n-1}(x;q)$. The  proofs of the other identities  are  similar and are omitted. Calculate the Askey-Wilson derivative of the two sides of \eqref{IM:def} and use that 
 $\mathcal{D}_q\mathcal{E}_q(x;y)=2\frac{q^{1/4}}{1-q}y \mathcal{E}_q(x;y)$.  This gives 
 \be \label{C:1}\sum_{n=0}^{\infty}\mathcal{D}_q\beta_n(x;q)w^n=2\frac{q^{1/4}}{1-q}w \sum_{n=0}^{\infty} \beta_n(x;q) w^n.\ee
 Therefore equating the coefficients of $w^n$ on the two sides of \eqref{C:1} gives  the required identity. 

\end{proof}
One can verify that 
\[\begin{split}\lim_{q\to 1}\mathcal{B}_n(x;q)(1-q)^{n-1}&=\frac{2^{2n-2}}{n!} B_{n}(\frac{x+1}{2}),\\
\lim_{q\to 1}\beta_n(x;q)(1-q)^{n-1}&=\frac{2^{2n-2}}{n!} B_{n}(\frac{x}{2}),
\end{split}
\]
and 
\[\begin{split}\lim_{q\to 1}\mathcal{E}_n(x;q)(1-q)^{n-1}&=\frac{2^{2n-1}}{n!} E_{n}(\frac{x+1}{2}),\\
\lim_{q\to 1}\widetilde{E}_n(x;q)(1-q)^{n-1}&=\frac{2^{2n}}{n!}B_{n}(\frac{x}{2}).
\end{split}
\]
\begin{rem} 
For $n\in\mathbb{N}$, set $T_n:=(\frac{2 q^{1/4}}{1-q})^n \dfrac{P_n}{[n]_q!}$,
where $P_n$ is any one of the polynomials $\mathcal{B}_n(x;q)$, $\beta_n(x;q)$, $\mathcal{E}_n(x;q)$, and $\widetilde{E}_n(x;q)$,  then we have 
\[\mathcal{D}_q T_n=[n]_q T_{n-1}.\]
That is, $(T_n)_n$ is a  sequence of $q$-Appell polynomials,  see~\cite{Al-Salam67, Tempesta}.
\end{rem}

Ismail~\cite[Eq. (14.6.4)]{Ismail-book}    defined  the action of the operator of  translation  by $y$,  $E_q^y$, on $H_n(x|q)$ as 
\be 
E_q^y H_n(x|q):=\sum_{m=0}^{n}\qbinom{n}{m} H_m(x|q) g_{n-m}(y)q^{\frac{m^2-n^2}{4}},
\ee
where 
\[g_n(x):=q^{n^2/4} (1+e^{2i\theta}) e^{-in \theta} (-q^{2-n}e^{2i\theta};q^2)_{n-1},\; x=\cos \theta.\]

He then extended  $E_q^y$ to the space of all polynomials by linearity. He also proves 
\[E_q^0=\mbox{identity}, \; g_n(0)=\delta_{n,0},\; \mbox{and}\; (E_q^y f)(x)=(E_q^x f)(y).\]

\begin{prop}\label{prop:B}   The $q$-Bernoulli polynomials  which are  defined in \eqref{S:def} and \eqref{IM:def} are related by the $q$-translation operator via the identity
\be \label{Eq:12}
E_q^{-\eta}\mathcal{B}_n(x;q)=\beta_n(x;q)
 \;(n\in\mathbb{N}_0).
\ee
Moreover, 
 \be\label{Eq:10} \mathcal{B}_n(q)=\beta_n(q)  \;  (n\in\mathbb{N}_0). \ee
 \end{prop}
 \begin{proof}
 Ismail proved that 
\[E_q^y \mathcal{E}_q(x;\alpha)=\mathcal{E}_q(x;\alpha)\mathcal{E}_q(y;\alpha),\]
see~\cite[Eq. (14.6.10)]{Ismail-book}. Therefore, acting by the operator $E_q^{-\eta} $ to the two sides of \eqref{S:def} yields 
\be\label{Eq:11}\sum_{n=0}^{\infty}E_q^{-\eta}\mathcal{B}_n(x;q) w^n=\frac{w\mathcal{E}_q(x;w)\mathcal{E}_q(-\eta;w)}{\mathcal{E}_q(\eta;w)-\mathcal{E}_q(-\eta;w)}
=\sum_{n=0}^{\infty}\beta_n(x;q) w^n.
\ee
Therefore, equating the coefficients of $w^n$ ($n\in\mathbb{N}_0$) on the two sides of \eqref{Eq:11} yields \eqref{Eq:12}. The  proof of  \eqref{Eq:10} follows by substituting with $x=-\eta $ on \eqref{S:def} and $x=0$ on \eqref{IM:def}. This yields 
\be\label{Eq:13} \begin{gathered}
\sum_{n=0}^{\infty}\mathcal{B}_n(-\eta;q) w^n=\dfrac{w\mathcal{E}_q(-\eta;w)}{\mathcal{E}_q(\eta;w)-\mathcal{E}_q(\eta;-w)}=\sum_{n=0}^{\infty}\beta_n(0;q) w^n.
\end{gathered}\ee
Equating the coefficients of $w^n$ on the two sides of \eqref{Eq:13}  yields $\mathcal{B}_n(-\eta;q)=\beta_n(0;q) $. This proves \eqref{Eq:10} and completes the proof of this  proposition.

 \end{proof}
\begin{rem}
Ismail and Rahman~\cite{Ismail-Rahman} introduced a $q$-analog of Bernoulli polynomials, $(B_n(x,y|q))_n$,   by 
\[
\begin{gathered}
-\sum_{n=0}^{\infty}B_n(x,y|q)\frac{t^n q^{n/4}}{(q;q)_n}
=\dfrac{(qt^2(1-q)^2;q^{2})_{\infty}}{(t^2(1-q)^2;q^2)_{\infty}}\mathcal{E}_q(x,-y;(1-q)t)\sum_{n=0}^{\infty}A_n(e^{i\theta},e^{i\phi})q^{n/4}(2t)^n\\
+t\dfrac{(qt^2(1-q)^2;q^)_{\infty}}{(t^2(1-q)^2;q^2)_{\infty}}\left[A_1(q^{1/2}e^{i\theta,e^{i\phi}})-A_1(e^{i\theta,e^{i\phi}})\right]  [\mathcal{E}_q(x,-y;(1-q)t)-\mathcal{E}_q(x,-y;-(1-q)t)].
\end{gathered}
\]
These  polynomials  rise through  their study for a fractional integral operator associated with the Askey-Wilson difference operator. The authors also pointed out that 
the particular choice,    $y=\eta=(q^{1/4}+q^{-1/4})/2$,   produces that $ B_n(x,\eta|q)$ is a constant multiplier of Suslov's $q$-Bernoulli polynomials  \eqref{S:def}.   Consequently,  from Proposition \ref{prop:B},  the polynomial 
$E_q^{-\eta}B_n(x,\eta|q)$ is a constant multiplier  of $\beta_n(x;q)\quad (n\in\mathbb{N})$.

\end{rem}
 \begin{prop}   The $q$-Euler  polynomials  defined in \eqref{S:ENdef} and \eqref{IM:Edef} are related by the $q$-translation operator via the identity
\be \label{Eq:14}
E_q^{-\eta}\mathcal{E}_n(x;q)=\frac{1}{2}\widetilde{E}_n(x;q),
\ee
 and 
 \be\label{Eq:15} \mathcal{E}_n(-\eta;q)=\frac{1}{2}\widetilde{E}_n(0;q), \ee
 where $n\in\mathbb{N}_0$. 
 \end{prop}
 \begin{proof}
 The proof is similar to the proof of Proposition~\ref{prop:B} and is omitted. 
 \end{proof}

\begin{prop} For $n\in\mathbb{N}_0$, 
\begin{eqnarray} \label{Eq:16}
\mathcal{B}_n(x;q)&=&\sum_{k=0}^{n}(-1)^{n-k}\beta_{n-k}(q)\dfrac{(-q^{1/4}e^{i\theta},-q^{1/4}e^{-i\theta};q^{1/2})_k}{(q;q)_k},\\\label{Eq:17}
\beta_n(x;q)&=&\sum_{k=0}^{n}\beta_{n-k}(q)\frac{q^{k^2/4}}{(q;q)_k}\rho_k(x),\\ \label{Eq:18}
H_n(x|q)&=&2q^{-n^2/4}(q;q)_n \sum_{k=0}^{n}\frac{q^{k^2+\frac{k}{2}}}{(q^{1/2};q^{1/2})_{2k+1}} \mathcal{B}_{n-2k}(x;q),
\end{eqnarray}

\end{prop}
\begin{proof}
To prove \eqref{Eq:16},  from the identity \cite[Eq. (14.1.2)]{Ismail-book}
\be
\mathcal{E}_q(x;w)=\mathcal{E}_q(\eta;w)\;{}_2\phi_1\left(-q^{1/4} e^{i\theta},-q^{1/4}e^{-i\theta};-q^{1/2};q^{1/2},w\right).
\ee
Thus,
\be
\sum_{n=0}^{\infty}\mathcal{B}_n(x;q)w^n=\dfrac{w\mathcal{E}_q(\eta;w)}{\mathcal{E}_q(\eta;w)-\mathcal{E}_q(\eta;-w)}{}_2\phi_1\left(-q^{1/4} e^{i\theta},-q^{1/4}e^{-i\theta};-q^{1/2};q^{1/2},w\right),
\ee
where $|w|<1$. 
Since 
\[\dfrac{w\mathcal{E}_q(\eta;w)}{\mathcal{E}_q(\eta;w)-\mathcal{E}_q(\eta;-w)}=\sum_{n=0}^{\infty}\beta_n(q)(-w)^n,\]
then 
\[\sum_{n=0}^{\infty}\mathcal{B}_n(x;q)w^n=\sum_{n=0}^{\infty}\beta_n (-w)^n,\]
and 
\be\label{Eq:19}\sum_{n=0}^{\infty}\mathcal{B}_n(x;q)w^n=\sum_{n=0}^{\infty}w^n\sum_{k=0}^{n}(-1)^{n-k}\beta_{n-k}(q)\dfrac{(-q^{1/4}e^{i\theta},-q^{1/4}e^{-i\theta};q)_k}{(q;q)_k}.\ee
Equating the coefficients of $w^n$ in \eqref{Eq:19} yields \eqref{Eq:16}. To prove \eqref{Eq:17}, note that from \eqref{IM:def},
\be\label{Eq:20}\begin{gathered}\sum_{n=0}^{\infty}\beta_n(x;q)w^n=\dfrac{w\mathcal{E}_q
(\eta;-w)}{\mathcal{E}_q(\eta;w)-\mathcal{E}_q(\eta;-w)}\mathcal{E}_q(x;w)\\
=\left(\sum_{n=0}^{\infty}\beta_n(q)w^n\right)\left( \sum_{n=0}^{\infty}\frac{q^{n^2/4}}{(q;q)_n} \rho_n(x)w^n\right)\\
=\sum_{n=0}^{\infty}w^n\sum_{k=0}^{n}\frac{q^{k^2/4}}{(q;q)_k} \rho_k(x)\beta_{n-k}(q).
\end{gathered}\ee
Therefore, equating  the coefficients of $w^n$ in the series on the sides of \eqref{Eq:20} yields \eqref{Eq:17}. 
Now, we prove \eqref{Eq:18}. From \eqref{Hermite-rep},  \eqref{S:def}, and 
\[\dfrac{(-w;q^{1/2})_{\infty}-(w;q^{1/2})_{\infty}}{w}=2\sum_{n=0}^{\infty}q^{\frac{n(2n+1)}{2}}\frac{w^{2n}}{(q^{1/2};q^{1/2})_{2n+1}},\]
we obtain
\be
2\left(\sum_{n=0}^{\infty}\mathcal{B}_n(x;q)w^n\right)\;\left(\sum_{n=0}^{\infty}\frac{q^{n^2+\frac{n}{2}}}{(q^{1/2};q^{1/2})_{2n+1}} w^{2n}\right)=\sum_{n=0}^{\infty}\frac{q^{n^2/4}}{(q;q)_n} H_n(x|q)w^n.
\ee
Consequently, 
\be\label{Eq:21}
2\sum_{n=0}^{\infty}w^n\sum_{k=0}^{[n/2]}\frac{q^{k^2+\frac{k}{2}}}{(q^{1/2};q^{1/2})_{2k+1}}\mathcal{B}_{n-2k}(x;q)=\sum_{n=0}^{\infty}\frac{q^{n^2/4}}{(q;q)_n} H_n(x|q)w^n.
\ee
Then, equating  the coefficients of $w^n$ in the series on the sides of \eqref{Eq:21} yields \eqref{Eq:18} and completes the proof of the proposition.  
\end{proof}

 In~\cite{Ismail-Mansour}, we   defined a $q$-analog of Bernoulli polynomials through the generating function 
\be
\dfrac{yE_q(xy)}{e_q(y/2)E_q(y/2)-1}=\sum_{n=0}^{\infty}B_n(x;q) \frac{y^n}{[n]!},
\ee
where $e_q(y):=\dfrac{1}{(y(1-q);q)_{\infty}},\quad E_q(y):=(-y(1-q);q)_{\infty}.$ The $q$-Bernoulli numbers  we defined  are the values of the Bernoulli polynomials at zero. I.e. 
\be\label{Eq:11-}
\dfrac{y}{e_q(y/2)E_q(y/2)-1}=\sum_{n=0}^{\infty}B_n(q) \frac{y^n}{[n]!}.
\ee
\begin{prop}
For $n\in\mathbb{N}_0$,
\[\mathcal{B}_n(q^2)=\beta_n(q^2)=B_n(q)\frac{2^{n-1}(1-q)}{(q;q)_{n}},\]
where $(\mathcal{B}_n(q))_n$ or $(\beta_n(q))_n$ are the $q$-Bernoulli numbers  in \eqref{S:def2} or \eqref{BND2} and $(B_n(q))_n$ are the $q$-Bernoulli number in \eqref{Eq:11-}. 
\end{prop}

\begin{proof}
The proof follows directly by comparing the  generating function in \eqref{S:def2} or \eqref{BND2} with the generating function in  \eqref{Eq:11-}. 
\end{proof}

In~ \cite{Ismail-Mansour}, we  proved that 
\[B_0(q)=1,\, B_1(q)=-\frac{1}{2},\,B_2(q)=\frac{q[2]}{2^2[3]},\; B_4(q)=-\frac{q^4}{2^4}\dfrac{(-q;q)_2[2]}{[3][5]}.\]
Moreover, $B_{2n+1}(q)=0$, $n=1,2,3,\ldots$, and $(-1)^{n-1}B_{2n}(q)>0$. Therefore,

\[\beta_0(q)=\frac{1-\sqrt{q}}{2},\, \beta_1(q)=\frac{-1}{2},\;\beta_2(q)= \frac{q^{1/2}}{2(1-q^{3/2})},
\]
\[\beta_{2n+1}(q)=0,\;(-1)^n\beta_{2n}(q)>0,\; n\in\mathbb{N}.\]

\begin{prop}
 The connection relations between the $q$-Bernoulli polynomials defined in \eqref{S:def} and \eqref{IM:def}  are given for 
 $n\in\mathbb{N}$ by 
\begin{eqnarray}
\label{F:1}
\beta_n(x;q)&=&\sum_{k=0}^{n}\frac{(-q^{-1/2};q)_k}{(q;q)_k} (-q^{1/2})^k \mathcal{B}_{n-k}(x;q),\\ \label{F:2}
\mathcal{B}_n(x;q)&=&\sum_{k=0}^{n}\frac{(-q^{1/2};q)_k}{(q;q)_k}  \beta_{n-k}(x;q).
\end{eqnarray}
\end{prop}

\begin{proof}
From \eqref{S:def2},  \eqref{BND2},  and
 \[\mathcal{E}_q(\eta;w)=\dfrac{(-w;q)_{\infty}}{(q^{1/2}w;q)_{\infty}}=\sum_{n=0}^{\infty}\dfrac{(-q^{-1/2};q)_n}{(q;q)_n} (q^{1/2}w)^n\;\left(q^{1/2}|w|<1\right),\]
we obtain 
\[\begin{split}
\sum_{n=0}^{\infty}\beta_n(x;q)  w^n&=\sum_{n=0}^{\infty}\mathcal{B}_n(x;q)  w^n \sum_{n=0}^{\infty}\dfrac{(-q^{-1/2};q)_n}{(q;q)_n} (-q^{1/2}w)^n\\
&=\sum_{n=0}^{\infty}w^n\sum_{k=0}^{n}\dfrac{(-q^{-1/2};q)_k}{(q;q)_k} (-q^{1/2})^k
\mathcal{B}_{n-k}(x;q).\end{split}\]
Equating the coefficients of $w^n$ yields \eqref{F:1}. Similarly, we can prove \eqref{F:2}.

\end{proof}

\begin{prop} For $n\in\mathbb{N}_0$
\[\begin{split}\mathcal{B}_n(-x;q)&=(-1)^n \mathcal{B}_n(x;q),\\
\beta_n(-x;q)&=(-1)^n\sum_{k=0}^{n}\frac{(-1;q^{1/2})_k}{(q^{1/2};q^{1/2})_k}\beta_{n-k}(x;q).
\end{split}
\]
\end{prop}
\begin{proof}
The proof follows from the generating functions \eqref{S:def} and \eqref{IM:def} and the facts that
\[\mathcal{E}_q(-x;y)=\mathcal{E}_q(x;-y),\]
\[\dfrac{\mathcal{E}_q(\eta;w)}{\mathcal{E}_q(\eta;-w)}=\dfrac{(-w;q^{1/2})_{\infty}}{(w;q^{1/2})_{\infty}}=\sum_{n=0}^{\infty}\dfrac{(-1;q^{1/2})_n}{(q^{1/2};q^{1/2})_n} w^n\quad (|w|<1).\]
\end{proof}

\section{A $q$ type Lidstone expansion in terms of $q$-Bernoulli polynomials }

 Ismail and Stanton~\cite[Theorem 3.4]{Ismail-Stanton-2003} proved the following theorem for entire functions of $q$-exponential growth of order less than $2 \ln q^{-1}$. However,  the theorem  is valid  if  f is an  entire function of $q$-exponential growth of order  $2 \ln q^{-1}$  and type less than $\frac{2\ln 2}{\ln q^{-1}}$.

\begin{thm}
Let $f$ be a function of $q$-exponential growth of order less than $ 2 \ln q^{-1}$ or of order $2 \ln q^{-1}$  and type less than $\frac{2\ln 2}{\ln q^{-1}}$. Then $f$ has the expansion 
 \be\label{Exp-2} f(x):=\sum_{k=0}^{\infty} {f_k} \rho_k(x),\;
 f_k:=\frac{q\frac{k^2-k}{4}(1-q)^k}{2^k (q;q)_k}\mathcal{D}_q^kf(0).\ee
 \end{thm}
 \begin{proof}
 The case when $f$ is a function of $q$-exponential growth of order less than $ 2 \ln q^{-1}$ is in \cite[Theorem 3.4]{Ismail-Stanton-2003}. Here, we assume  that  $f$ is a function of $q$-exponential growth of order $ 2 \ln q^{-1}$ and type $\alpha$. The proof is similar to the proof in \cite[Theorem 3.4]{Ismail-Stanton-2003} , the critical step is the  step of proving that $I_m \to 0 $ as $m\to\infty$, where  $I_m$ satisfies the inequality 
 \[|I_m|\leq B M(r_m;f) q^{\frac{(m+\delta)^2}{4}},\;
r_m:=\dfrac{q^{-\frac{m+\delta}{2}}-q^{\frac{m+\delta}{2}}}{2},\]
where $m$ is an even positive integer, $\delta$ and $B$ are positive  constants. 
Consequently, 
\[\ln r_m=-\frac{m+\delta}{2} \ln q -\ln 2+\ln(1-q^{m-\delta}),\]
\[\frac{\ln^2 r_m}{\ln q}=\frac{(m+\delta)^2}{4}\ln q+(m+\delta)\ln 2+\mathcal{O}(1),\]
as $m\to \infty$. 
Since $M(r_m;f)\leq C r_m^{\alpha} e^{-\frac{\ln^2 r_m}{\ln q}}$, $C$ is a positive constant,  then 
\[\ln |I_m|\leq  -(m+\delta) (\ln 2+\frac{\alpha}{2}\ln q)+\mathcal{O}(1),\]
as $m\to\infty$. Therefore, if $(\ln 2+\frac{\alpha}{2}\ln q)>0$, i.e. $\alpha<\frac{2\ln 2}{\ln q^{-1}}$,  $\lim_{m\to\infty}I_m=0$ as required.

 \end{proof}
\begin{thm}\label{thm:1}
Let $f$ be a function of $q$-exponential growth of order less than $ 2 \ln q^{-1}$ or of order $2 \ln q^{-1}$  and type less than $\frac{2\ln 2}{\ln q^{-1}}$.  Then $f$ has the expansion 
 \be\label{Exp-1}f(x):=\sum_{k=0}^{\infty} {f_k} \rho_k(x),\;
 f_k:=\frac{q\frac{k^2-k}{4}(1-q)^k}{2^k (q;q)_k}\mathcal{D}_q^kf(0).\ee
 Let
\[F(y):=\sum_{k=0}^{\infty} \frac{f_k}{\psi_k}\frac{1}{y^{k+1}},\;\psi_k:=\frac{q^{k^2/4}}{(q;q)_k}.\]
If  $\displaystyle\limsup_{n\to\infty}\sqrt[n]{\left|\mathcal{D}_q^nf(0)\right|}<\frac{2q^{1/4}}{1-q}$, 
then
$\displaystyle \tau:=\limsup_{n\to\infty}\sqrt[n]{\left| \frac{f_n}{\psi_n}\right|}<1$, 
and  
\[f(x):=\frac{1}{2\pi\,i} \int_{\Gamma}\mathcal{E}_{q}(x;y) F(y)\,dy,\]
where $\Gamma$ is a  circle centered  at zero  and of radius $r$, $\tau<r<1$.

\end{thm}

\begin{proof}

The series  $F(y)$ is uniformly convergent for   $\Gamma:=\{y\in\mathbb{C}:|y|=r\}$  since $\tau<r<1$. Consequently,
\be
\frac{1}{2\pi\,i}\int_{\Gamma}\mathcal{E}_q(x;y)F(y)\,dy=\sum_{k=0}^{\infty}\frac{f_k}{\psi_k} \frac{1}{2\pi\,i}\int_{\Gamma}\frac{\mathcal{E}_q(x;y)}{y^{k+1}}\,dy.
\ee
But from the Cauchy Residue Theorem 
\be 
\frac{k!}{2\pi\,i}\int_{\Gamma}\frac{\mathcal{E}_q(x;y)}{y^{k+1}}\,dy=\frac{d^k}{dy^k}\mathcal{E}_q(x;y)|_{y=0}.
\ee
From \eqref{rho-rep}, we obtain  
\[\mathcal{E}_q(x;y)=\sum_{k=0}^{\infty}\psi_k \rho_k(x)y^k.\]
Consequently 
\be 
\frac{d^k}{dy^k}\mathcal{E}_q(x;y)|_{y=0}=k! \psi_k \rho_k(x),
\ee
and 
\[\frac{1}{2\pi\,i}\int_{\Gamma}\mathcal{E}_q(x;y)F(y)\,dy=\sum_{k=0}^{\infty}f_k \rho_k(x)=f(x).\]
This completes the proof of the theorem.
\end{proof}

\begin{thm} \label{thm:M0} Let $\gamma:=\frac{q^{1/4}(1-q)}{2}$. Then 
\be\label{Eq:4}
\mathcal{E}_q(x;y)=-\sum_{k=0}^{\infty} \gamma^{2k} B_k(x) y^{2k}+\mathcal{E}_q(\eta;y) \sum_{k=0}^{\infty} \gamma^{2k}A_k(x) y^{2k},
\ee
where $\{A_k(x)\}$ and $\{B_k(x)\}$ are the polynomials generated by 

\be\label{Eq:1-2}
\begin{split}
\sum_{k=0}^{\infty}\gamma^{2k} A_k(x) y^{2k}&=\frac{\mathcal{E}_q(x;y)-\mathcal{E}_q(x;-y)}{\mathcal{E}_q(\eta;y)-\mathcal{E}_q(\eta;-y)}\\
&=(qy^2;q^2)_{\infty}\frac{\mathcal{E}_q(x;y)-\mathcal{E}_q(x;-y)}{(-y;q^{1/2})_{\infty}-(y;q^{1/2})_{\infty}},
\end{split}
\ee
\be\begin{split} \label{Eq:1}
\sum_{k=0}^{\infty}\gamma^{2k} B_k(x) y^{2k}&=\frac{\mathcal{E}_q(\eta;-y)\mathcal{E}_q(x;y)-\mathcal{E}_q(\eta;y)\mathcal{E}_q(x;-y)}{\mathcal{E}_q(\eta;y)-\mathcal{E}_q(\eta;-y)}\\
&=\frac{(y;q^{1/2})_{\infty}\mathcal{E}_q(x;y)-(-y;q^{1/2})_{\infty}\mathcal{E}_q(x;-y)}{(y;q^{1/2})_{\infty}-(-y;q^{1/2})_{\infty}}.
\end{split}\ee

The series in \eqref{Eq:1-2}  and \eqref{Eq:1} converge uniformly on any compact subsets of $|y|<w_1$, where $w_1$ is the   smallest positive zeros of  $S_q(\eta;y)$.
Moreover, 
\be\label{57}
A_k(x)=2\gamma^{-2k} \mathcal{B}_{2k+1}(x;q),\quad B_k(x)=2\gamma^{-2k}\beta_{2k+1}(x;q),
\ee
where $(\mathcal{B}_k(x;q))_k$ and $(\beta_k(x;q))_k$ are the $q$-Bernoulli polynomials defined by the generating functions in \eqref{S:def} and \eqref{IM:def}, respectively. 
\end{thm}
\begin{proof}
The proof follows from the identity 
\[\mathcal{E}_q(x;y)=-\frac{\mathcal{E}_q(\eta;-y)\mathcal{E}_q(x;y)-\mathcal{E}_q(\eta;y)\mathcal{E}_q(x;-y)}{\mathcal{E}_q(\eta;y)-\mathcal{E}_q(\eta;-y)}+\mathcal{E}_q(\eta;y)\frac{\mathcal{E}_q(x;y)-\mathcal{E}_q(x;-y)}{\mathcal{E}_q(\eta;y)-\mathcal{E}_q(\eta;-y)},\]
and the generating functions in \eqref{Eq:1-2} and \eqref{Eq:1}. 
\end{proof}

\begin{thm}\label{Thm:M}
Let  $\eta:=\frac{q^{1/4}+q^{-1/4}}{2}$, and $w_1$ be the smallest positive zero of $S_q(\eta,w)$.   Assume  that  $f$ is  a function of $q$-exponential growth of order less than $2\ln q^{-1}$ or of order $2\ln q^{-1}$ and of type less than $\frac{2\ln 2}{\ln q^{-1}}$.  
 Then $f$ has the expansion 
 \be\label{Exp}f(x):=\sum_{k=0}^{\infty} {f_k} \rho_k(x),\;
 f_k:=\frac{q\frac{k^2-k}{4}(1-q)^k}{2^k (q;q)_k}\mathcal{D}_q^kf(0).\ee
If \be \label{cond}\limsup_{n\to\infty}\sqrt[n]{|\mathcal{D}_q^nf(0)|q^{-n/4}}=\tau\;, 0<\tau <\min( 1,w_1),\ee
then 
\be \label{expansion}f(x)=2\sum_{k=0}^{\infty}\gamma^{-2k}  \mathcal{D}_q^{2k}f(\eta)\mathcal{B}_{2k+1}(x;q)-2\sum_{k=0}^{\infty} \gamma^{-2k}\mathcal{D}_q^{2k}f(0)\beta_{2k+1}(x;q) , \quad x\in\mathbb{C},\ee
where $\mathcal{B}_k(x;q)$ and $\beta_k(x;q)$ are the $q$-Bernoulli polynomials defined by the generating functions in \eqref{S:def} and \eqref{IM:def}, respectively. 
\end{thm}

\begin{proof}
From Theorem \ref{thm:1},
\[f(x)=\frac{1}{2\pi\,i}\int_{\Gamma} \mathcal{E}_q(x;y) F(y)\,dy\quad (x\in\mathbb{C}),\]
where $\Gamma$ is a circle centered at the origin of radius $r$, $\tau<r<1$. Then
\be 
\mathcal{D}_q^{2k} f(x)=\frac{1}{2\pi\,i}\int_{\Gamma}\gamma^{2k} y^{2k} \mathcal{E}_q(x;y) F(y)\,dy.
\ee
Hence 
\[\begin{split}\mathcal{D}_q^{2k} f(0)&=\frac{1}{2\pi\,i}\int_{\Gamma}\gamma^{2k} y^{2k}  F(y)\,dy,\\
\mathcal{D}_q^{2k} f(\eta)&=\frac{1}{2\pi\,i}\int_{\Gamma}\gamma^{2k} \mathcal{E}_q(\eta;y)y^{2k}  F(y)\,dy,
\end{split}
\]
and for all $n\in\mathbb{N}$, we have 
\be\label{Eq:n}\begin{gathered}-\sum_{k=0}^{n}\mathcal{D}_q^{2k} f(0) B_k(x)+\sum_{k=0}^{n}\mathcal{D}_q^{2k} f(\eta) A_k(x)=\\
\frac{1}{2\pi\,i}\int_{\Gamma}\left(- \sum_{k=0}^{n}\gamma^{2k} B_k(x)y^{2k}+\mathcal{E}_q(\eta;y)\sum_{k=0}^{n}
\gamma^{2k}A_k(x)y^{2k}\right) F(y)\,dy.
\end{gathered} \ee
From the generating functions \eqref{Eq:1} and \eqref{Eq:1-2} define the polynomials  $(A_k)$ and $(B_k)$, the series\[ \sum_{k=0}^{\infty}\gamma^{2k} A_k(x)y^{2k}\;\mbox{ and}\; \sum_{k=0}^{n}\gamma^{2k}B_k(x)y^{2k} \] are uniformly convergent for $|y|<w_1$, where $w_1$ is the smallest positive zero of $S_q(\eta;w)$. Moreover, for $|y|<\min (1,w_1)$,
\[- \sum_{k=0}^{n}\gamma^{2k}B_k(x)y^{2k} +\mathcal{E}_q(x;\eta)
\sum_{k=0}^{\infty}\gamma^{2k} A_k(x)y^{2k}\] converges uniformly to $\mathcal{E}_q(x;y)$.  So,  we can take the limit as $n\to \infty$ on \eqref{Eq:n} to obtain 
\be \begin{gathered}\sum_{k=0}^{\infty}\left(\mathcal{D}_q^{2k} f(\eta) A_k(x)-\mathcal{D}_q^{2k} f(0) B_k(x)\right)=
\frac{1}{2\pi\,i}\int_{\Gamma }\mathcal{E}_q(x;y)F(y)\,dy=f(x).
\end{gathered} \ee
\end{proof}

\vskip .5 cm 

\begin{exa}
This  example shows 
that the condition on the constant $\tau$ appears in \eqref{cond} in  Theorem \ref{Thm:M} is essential. If we take $f(x)=S_q(x;w_1)$ then $\tau=w_1$ and 
\[\mathcal{D}_q^{2k} f(0)=\mathcal{D}_q^{2k} f(\eta)=0\quad (k\in\mathbb{N}_0).\] Therefore,   the function can not be expended as in \eqref{expansion}.

\end{exa}

\vskip .5 cm

\begin{thm}
The polynomials $\left(A_k(x)\right)_k$ and $\left(B_k(x)\right)_k$defined  in \eqref{Eq:1} and \eqref{Eq:1-2}.
 satisfy the following identities for $k\in\mathbb{N}$. 

\begin{itemize} 
\item[(I) ] $E_{q}^{-\eta}B_k(x)=A_k(x)$.
\item [(II) ] $A_0(x)\equiv 1$,  $\mathcal{D}_q^2 A_k(x)=A_{k-1}(x)$.

 \item[(III)] $B_0(x)\equiv 1$,  $\mathcal{D}_q^2 B_k(x)=B_{k-1}(x)$.
 \item[(IV)] $A_k(0)=B_k(0)=0$.
 \item[(V)]$A_k(1)=B_k(1)=0$.
 \end{itemize}
 
 \end{thm}
\begin{proof}
The proof  of (I)  follows from \eqref{Eq:12} and \eqref{57}. 
The proofs of (II) and (III) follow by acting by  the Askey-Wilson operator  of the two sides of the generating function of $(A_k(x))_k$ and $(B_k(x))_k$
and using that 
\[\mathcal{D}_q \mathcal{E}_q(x;y)=\frac{2q^{1/4} y}{1-q} \mathcal{E}_q(x;y).\]
The proofs of (IV) and (V) follow by substituting $x=0$ and $x=1$ on the generating functions identities of $(A_k(x))$ and $(B_k(x))$ and using that 
$\mathcal{E}_q(0;\alpha)=1$.
\end{proof}
\begin{exa}\label{Ex:4}
Consider the function $\phi_n(x;a)=(ae^{i\theta},ae^{-i\theta};q)_{n}$, $n\in\mathbb{N}$, and $x=\cos\theta$.
Since
\[\mathcal{D}_q\phi_n(x)=-2aq^{-1/2}[n]_q \phi_{n-1}(x;aq^{1/2}),\]
then 
\[\mathcal{D}_q^{2k}\phi_{2n}(x)=(2a)^{2k} q^{k^2-\frac{k}{2}}\dfrac{[2n]_q!}{[2n-2k]_q!} \phi_{2n-2k}(x;aq^{k})\quad (k=0,1,\ldots,n),\]
and zero otherwise. 
Hence, for $k=0,1,\ldots,n$,
\[\mathcal{D}_q^{2k}\phi_{2n}(0)=(2a)^{2k} q^{k^2-\frac{k}{2}}\dfrac{[2n]_q!}{[2n-2k]_q!} (-a^2 q^{2k};q^2)_{2n-2k},\]
\[\mathcal{D}_q^{2k}\phi_{2n}(\eta)=(2a)^{2k} q^{k^2-\frac{k}{2}}\dfrac{[2n]_q!}{[2n-2k]_q!} (a q^{k}q^{1/4}, aq^{k}q^{-1/4};q)_{2n-2k},\]
and 
\[\begin{gathered}(ae^{i\theta}, ae^{-i\theta};q)_{2n}=2\sum_{k=0}^{n}\gamma^{-2k}(2a)^{2k} q^{k^2-\frac{k}{2}}\dfrac{[2n]_q!}{[2n-2k]_q!} (a q^{k}q^{1/4}, aq^{k}q^{-1/4};q)_{2n-2k}\mathcal{B}_{2k+1}(x;q)\\
-2\sum_{k=0}^{n}\gamma^{-2k}(2a)^{2k} q^{k^2-\frac{k}{2}}\dfrac{[2n]_q!}{[2n-2k]_q!} (-a^2 q^{2k};q^2)_{2n-2k}
\mathbb{\beta}_{2k+1}(x;q).
\end{gathered}\]

\end{exa}
\section{A $q$-Lidstone expansion in terms of $q$-Euler polynomials}
Whittaker in \cite{Whittaker} expanded the entire functions that satisfy the condition
\[\limsup_{r\to\infty}\frac{\ln M(r;f)}{r}<\frac{\pi}{2}\]as 
\[f(z)=\sum_{n=0}^{\infty} f^{(2n)}(1) M_n(z)-\sum_{n=0}^{\infty} f^{(2n+1)}(0)M_{n+1}'(1-z), \]
where $(M_n)_n$ is the sequence of polynomials generated by 
\[\frac{e^{zt}+e^{-zt}}{e^{t}+e^{-t}}=\sum_{n=0}^{\infty}M_n(z) t^{2n}. \]
It turned out that $M_n(z)=\frac{2^{2n}}{2n!}E_{2n}(\frac{z+1}{2})$, where $(E_n(z))_n$ is the sequence of Euler's polynomials defined as in \eqref{B-E}. In this section, we extend Whittaker's result to the Askey-Wilson $q$-difference operator.

\begin{thm} \label{thm:E0} Let $\gamma:=\frac{q^{1/4}(1-q)}{2}$. Then 
\be\label{Eq:E4}
\mathcal{E}_q(x;y)=\sum_{k=0}^{\infty} \gamma^{2k} M_k(x) y^{2k}+\mathcal{E}_q(\eta;y) \sum_{k=0}^{\infty} \gamma^{2k+1}\widetilde{M}_k(x) y^{2k+1},
\ee
where $\{M_k(x)\}$ and $\{\widetilde{M}_k(x)\}$ are the polynomials generated by

\be\begin{split} \label{E:Eq2}
\sum_{k=0}^{\infty}\gamma^{2k+1} M_k(x) y^{2k+1}&=\frac{\mathcal{E}_q(\eta;-y)\mathcal{E}_q(x;y)-\mathcal{E}_q(\eta;y)\mathcal{E}_q(x;-y)}{\mathcal{E}_q(\eta;y)+\mathcal{E}_q(\eta;-y)}\\
&=\frac{(y;q^{1/2})_{\infty}\mathcal{E}_q(x;y)-(-y;q^{1/2})_{\infty}\mathcal{E}_q(x;-y)}{(y;q^{1/2})_{\infty}+(-y;q^{1/2})_{\infty}},
\end{split}\ee
\be\label{E:Eq1}
\begin{split}
\sum_{k=0}^{\infty}\gamma^{2k} \widetilde{M}_k(x) y^{2k}&=\frac{\mathcal{E}_q(x;y)+\mathcal{E}_q(x;-y)}{\mathcal{E}_q(\eta;y)+\mathcal{E}_q(\eta;-y)}\\
&=(qy^2;q^2)_{\infty}\frac{\mathcal{E}_q(x;y)+\mathcal{E}_q(x;-y)}{(-y;q^{1/2})_{\infty}+(y;q^{1/2})_{\infty}}.
\end{split}
\ee
The series in \eqref{Eq:E4} converges uniformly on any compact subsets of $|y|<\widetilde{w}_1$, where $\widetilde{w}_1$ is  smallest positive zeros of  $C_q(\eta;y)$.

\end{thm}
\begin{proof}
The proof follows from the identity 
\[\mathcal{E}_q(x;y)=\frac{\mathcal{E}_q(\eta;-y)\mathcal{E}_q(x;y)-\mathcal{E}_q(\eta;y)\mathcal{E}_q(x;-y)}{\mathcal{E}_q(\eta;y)+\mathcal{E}_q(\eta;-y)}+\mathcal{E}_q(\eta;y)\frac{\mathcal{E}_q(x;y)+\mathcal{E}_q(x;-y)}{\mathcal{E}_q(\eta;y)+\mathcal{E}_q(\eta;-y)},\]
and the generating functions in \eqref{E:Eq2} and \eqref{E:Eq1}. 

\end{proof}
\begin{prop}
If $ (M_k)_k$ and $(\widetilde{M}_k)$ are the polynomials defined by the generating functions  in \eqref{E:Eq2} and \eqref{E:Eq1}, respectively, then 
\[M_k(x)=\gamma^{-2k-1}\widetilde{E}_{2k+1}(x;q),\quad  \widetilde{M}_k(x)=2\gamma^{-2k}\mathcal{E}_{2k}(x;q),\]
where $\mathcal{E}_{n}(x;q)$ and $\widetilde{E}_n(x;q)$ are the $q$-Euler polynomials defined in \eqref{S:ENdef}, and \eqref{IM:Edef}, respectively. 
\end{prop}
\begin{proof}
The proof is straightforward and is omitted. 
\end{proof}

\begin{thm}\label{Thm:M2}
Let  $\eta:=\frac{q^{1/4}+q^{-1/4}}{2}$.   Assume  that  $f$ is  a function of $q$-exponential growth of order less than $2\ln q^{-1}$ or of order  $2\ln q^{-1}$ and type less than $\frac{2\ln 2}{\ln q^{-1}}$.
 Then $f$ has the expansion 
 \be\label{Exp2}f(x):=\sum_{k=0}^{\infty} {f_k} \rho_k(x),\;
 f_k:=\frac{q\frac{k^2-k}{4}(1-q)^k}{2^k (q;q)_k}\mathcal{D}_q^kf(0).\ee
If \be \label{cond2}\limsup_{n\to\infty}\sqrt[n]{|\mathcal{D}_q^nf(0)|q^{-n/4}}=\tau\;, 0<\tau <\min( 1,\widetilde{w}_1),\ee
where $\widetilde{w}_1$ is the smallest positive zero of $C_q(\eta;w)$, 
then 
\be \label{expansion2}f(x)=\sum_{k=0}^{\infty}\gamma^{-2k-1}\widetilde{E}_{2k+1}(x;q) \mathcal{D}_q^{2k+1}f(0)+
2\sum_{k=0}^{\infty}\gamma^{-2k}\mathcal{E}_{2k}(x;q)\mathcal{D}_q^{2k}f(\eta), \quad x\in\mathbb{C},\ee
 where $\mathcal{E}_{n}(x;q)$ and $\widetilde{E}_n(x;q)$ are the $q$-Euler polynomials defined in \eqref{S:ENdef}, and \eqref{IM:Edef}, respectively. 

\end{thm}
\begin{proof}
From Theorem \ref{thm:1},
\[f(x)=\frac{1}{2\pi\,i}\int_{\Gamma} \mathcal{E}_q(x;y) F(y)\,dy\quad (x\in\mathbb{C}),\]
where $\Gamma$ is a circle centered at the origin of radius $r$, $\tau<r<1$. Then
\be 
\mathcal{D}_q^{k} f(x)=\frac{1}{2\pi\,i}\int_{\Gamma}\gamma^{k} y^{k} \mathcal{E}_q(x;y) F(y)\,dy,
\ee

Hence 
\[\mathcal{D}_q^{2k+1} f(0)=\frac{1}{2\pi\,i}\int_{\Gamma}\gamma^{2k+1} y^{2k+1}  F(y)\,dy,\]
\[\mathcal{D}_q^{2k} f(\eta)=\frac{1}{2\pi\,i}\int_{\Gamma}\gamma^{2k} \mathcal{E}_q(\eta;y)y^{2k}  F(y)\,dy,\]
and for all $n\in\mathbb{N}$, we have 
\be\label{E:Eq:n} \begin{gathered}\sum_{k=0}^{n}\mathcal{D}_q^{2k+1} f(0) M_k(x)+\sum_{k=0}^{n}\mathcal{D}_q^{2k} f(\eta) \widetilde{M}_k(x)=\\
\frac{1}{2\pi\,i}\int_{\Gamma}\left(\sum_{k=0}^{n}\gamma^{2k+1} M_k(x)y^{2k+1}+\mathcal{E}_q(\eta;y)\sum_{k=0}^{n}
\gamma^{2k}\widetilde{M}_k(x)y^{2k}\right) F(y)\,dy.
\end{gathered} \ee
From the generating functions \eqref{E:Eq1} and \eqref{E:Eq2} define the polynomials  $(M_k)_k$ and $(\widetilde{M}_k)_k$, the series\[ \sum_{k=0}^{\infty}\gamma^{2k} M_k(x)y^{2k}\;\mbox{ and}\; \sum_{k=0}^{n}\gamma^{2k}\widetilde{M}_k(x)y^{2k} \] are uniformly convergent for $|y|<\widetilde{w}_1$, where $\widetilde{w}_1$ is the smallest positive zero of $C_q(\eta;w)$. Moreover, for $|y|<\min (1,\widetilde{w}_1)$, the series are uniformly convergent and we can take the limit as $n\to \infty$ on \eqref{E:Eq:n} to obtain 
\be \begin{gathered}\sum_{k=0}^{\infty}\left(\mathcal{D}_q^{2k+1} f(0) M_k(x)+(\mathcal{D}_q^{2k} f(\eta) \widetilde{M}_k(x)\right)=
\frac{1}{2\pi\,i}\int_{\Gamma }\mathcal{E}_q(x;y)F(y)\,dy=f(x).
\end{gathered} \ee
\end{proof}

\begin{exa}
In this example,  we show that the condition on the constant $\tau$ appears in \eqref{cond2} in  Theorem \ref{Thm:M2} is essential. If we take $f(x)=C_q(x;\widetilde{w}_1)$ then $\tau=\widetilde{w}_1$ and 
$\mathcal{D}_q^{2k+1} f(0)=\mathcal{D}_q^{2k} f(\eta)=0$ and  the function can not be expended as in \eqref{expansion2}.

\end{exa}

\begin{exa}
Consider the function $\phi_n(x;a)=(ae^{i\theta},ae^{-i\theta};q)_{n}$, $n\in\mathbb{N}$, and $x=\cos\theta$. similar to Example \ref{Ex:4}, we can prove 
\[\mathcal{D}_q^{2k+1}\phi_{2n}(0)=-(2a)^{2k+1} q^{k^2+\frac{k}{2}}\dfrac{[2n]_q!}{[2n-2k-1]_q!} (-a^2 q^{2k+1};q^2)_{2n-2k-1}\quad (k=0,1,\ldots,n-1),\]
\[\mathcal{D}_q^{2k}\phi_{2n}(\eta)=(2a)^{2k} q^{k^2-\frac{k}{2}}\dfrac{[2n]_q!}{[2n-2k]_q!} (a q^{k}q^{1/4}, aq^{k}q^{-1/4};q)_{2n-2k}\quad (k=0,1,2,\ldots,n),\]
and  they are equal to zero otherwise.  Consequently, 
\[\begin{gathered}(ae^{i\theta}, ae^{-i\theta};q)_{2n}=2\sum_{k=0}^{n}\gamma^{-2k}(2a)^{2k} q^{k^2-\frac{k}{2}}\dfrac{[2n]_q!}{[2n-2k]_q!} (a q^{k}q^{1/4}, aq^{k}q^{-1/4};q)_{2n-2k}\mathcal{E}_{2k}(x;q)\\
-\sum_{k=0}^{n-1}\gamma^{-2k-1}(2a)^{2k+1} q^{k^2+\frac{k}{2}}\dfrac{[2n]_q!}{[2n-2k-1]_q!} (-a^2 q^{2k+1};q^2)_{2n-2k-1}
\widetilde{E}_{2k+1}(x;q).
\end{gathered}\]

\end{exa}
\section{A $q$-analog of Guichard- Whittaker's Interpolation problem}

Whittaker, in his book~\cite{Whittaker},  raised the question :

Given an entire  function $f$, is there an entire function $g$ such that 
\be \label{G}g(z+1)-g(z)=f(z)?\ee
The  answer was yes,  if 
\[f(z)=\sum_{n=0}^{\infty}a_n z^n,\;\mbox{ then}\; g(z)=\sum_{n=0}^{\infty}\frac{a_n}{n+1}B_{n+1}(z),\]
where $B_n(z)$ is the classical Bernoulli number of order $n$. Guichard,  see~ \cite{Guichard},   was the first one who proved the result, but  Whittaker gave a more precise result in \cite{Whittaker}. 

In this section,  we introduce some essential concepts then we raise a similar question. Let  $p$ be  a positive number, define 
 a transformation operator $T_p^1$ acting on a polynomial $x^n$ as 
\be \label{T:def}
T_p^1 x^n:=T_n(x):=\sum_{k=0}^{n}\pbinom{n}{k} x^{n-k} \delta_k(p), \quad \lim_{p\to 1} \delta_k(p)=1,
\ee
and acting on power series $h(x):=\sum_{n=0}^{\infty}a_n x^n$ as a linear operator. I.e. 
\[ h(x \dotplus 1)=T_p^1 h(x)=\sum_{n=0}^{\infty}a_n T_n(x). \]
We define the function $E(t;p)$ by
$E(t;q):=\sum_{n=0}^{\infty}  \frac{t^n}{[n]_p
!}$ . Therefore,  
\[E(t;p)=\left\{\begin{array}{cc}e_q(t),&q=p<1,\\ \exp(t),&p=1,\\
E_q(t),&p=\frac{1}{q}; q<1.\end{array}\right.\]
\[ T_1^1 x^n=(x+1)^n\quad \mbox{ and }\quad h(x\dotplus 1)=h(x+1).\]

\vskip .5 cm 

The problem is as the following:  Given an entire function $f$,  is there an entire  function $g$ such that 
\be\label{FEP} g(z\dotplus 1)-g(z)=f(z)?\ee

To solve the problem, we try to find  a  sequence of polynomials $\phi_n(z)$ that satisfy the conditions 
\be \label{Eq1:phn}\phi_n(z\dotplus 1)-\phi_n(z)=[n]_p z^{n-1}, \quad [n]_p=\frac{1-p^n}{1-p} ,\ee
\be\label{Eq2:phn} \phi_0=1,\;  \mbox{ and}\;  \phi_n(0)=0 \quad ( n\in\mathbb{N}).\ee
    From \eqref{Eq1:phn}, 
\be\label{Eq3:phn} D_{p,z}^k[\phi_n(z\dotplus 1)-\phi_n(z)]_{z=0}=[n]_p! \delta_{k,n-1}, \ee
where $\delta_{k,n-1}$ is the Kronecker's delta.
\begin{prop}\label{I:prop-1}
Let $f$ be a polynomial of order $n$.  Assume that $(\phi_k)_k$ is a sequence of polynomials satisfying \eqref{Eq1:phn} and \eqref{Eq2:phn}. Then 
\be
f(z)=f(0)+\sum_{k=1}^{n}\frac{D_p^{k-1}f(0\dotplus 1)-D_p^{k-1}f(0)}{[k]_p!} \phi_k(z).
\ee

\end{prop} 
\begin{proof}
Assume that $f(z)=\sum_{k=0}^{n}b_k\phi_k(z)$. Then $f(0)=b_0$ because $\phi_k(0)=0$ for $k=1,2,\ldots,m$. 
Then 
\[D_p^m f(0\dotplus 1)-D_p^m f(0)=\sum_{k=1}^{m} b_kD_{p}^m[\phi_k(z\dotplus 1)-\phi_k(z)]_{z=0}=\sum_{k=1}^{m}\delta_{m,k-1}[k]_p!=b_{m+1}[m+1]_p!.  \]
Hence $b_m=\dfrac{D_p^{m-1} f(0\dotplus 1)-D_p^{m-1} f(0)}{[m]_p!}$, $m=1,2,\ldots,n$.
\end{proof}
\begin{lem}
If $f(z)$ is an entire  function then 
\[ D_{p,z} f(z\dotplus 1)=(D_{p,z}f)(z\dotplus 1).\]
\end{lem}
\begin{proof}
Assume that $f(z)=\sum_{k=0}^{\infty}a_k z^k$. Then 
$D_pf(z)=\sum_{k=1}^{\infty}a_k[k]_pz^{k-1}$. Consequently, 
\[(D_p f)(z\dotplus 1)=\sum_{k=1}^{\infty}a_k [k]_p\sum_{j=0}^{k-1}\pbinom{k-1}{j}\delta_j(p) z^{k-1-j}.\]
On the other hand, 
\[D_p f(z\dotplus  1)=D_{p,z}\sum_{k=0}^{\infty}a_k \sum_{j=0}^{k}\pbinom{k}{j } \delta_j(p)z^{k-j}=\sum_{k=1}^{\infty}a_k [k]_p\sum_{j=0}^{k-1}\pbinom{k-1}{j}\delta_j(p) z^{k-1-j}.\]
This proves the required result.
\end{proof}
\begin{prop}
Let $(\phi_n)_n$ be the set of polynomials defined by \eqref{Eq1:phn}-\eqref{Eq3:phn}. Then 
\[D_p\phi_n(x)=D_p\phi_n(0)+[n]_p \phi_{n-1}(x).\]

\end{prop}

\begin{proof}
We apply Proposition \ref{I:prop-1} to the function $\phi_n(x)$ ($n\in\mathbb{N}$). Then 
\[\begin{split}D_p\phi_n(x)&=D_p\phi_n(0)+\sum_{k=1}^{\infty} \dfrac{D_p^k\phi_n(0\dotplus 1)-D_p^k \phi_n(0)}{[k]_p!}\phi_k(x)\\
&=D_p\phi_n(0)+\sum_{k=1}^{\infty}\dfrac{[n]_p!}{[k]_p!}\delta_{k,n-1}\phi_k(x)\\
&=D_p\phi_n(0)+[n]_p\phi_{n-1}(x).
\end{split}\]

\end{proof}

Now we define  the polynomials  $(B_{p,n}(z))_n$ in terms of the polynomials $(\phi_n(z))_n$  by 
\[B_{p,n}(z):=\Phi_n(z)+\dfrac{D_p\phi_n(0)}{[n]_p},\]
We can prove that \[D_pB_{p,n}(z)=[n]_pB_{p,n-1}(z),\]

\be \label{Eq 4}D_{p,z}^k[B_{p,n}(z\dotplus 1)-B_{p,n}(z)]_{z=0}=[n]_p! \delta_{k,n-1} \quad (k=0,1,2,\ldots,n),\ee
\[B_{p,n}(z)=\sum_{k=0}^{n} \pbinom{n}{k} B_{n-k}(0)z^k,\quad B_{k}:=B_{p,k}(0).\]

\begin{thm}
The generating function of the polynomials $B_{p,n}(z)$ is given by 
\be 
\sum_{n=0}^{\infty}B_{p,n}(z) \frac{t^n}{[n]_p!}=  \frac{tE(tz;p)}{E(t(z\dotplus 1);p)|_{z=0}-1},
\ee
where $E(tz;p)=\sum_{k=0}^{\infty} \dfrac{(tz)^k}{[k]_p!}$ $(p\neq 1)$,
\[ E(t(z\dotplus 1);p)|_{z=0}=\sum_{k=0}^{\infty}\frac{t^k}{[k]_p!}\delta_k,\]
and  $(\delta_k)_k$ is the sequence defined in \eqref{T:def}. 
 
\end{thm}
\begin{proof}
Simple calculation shows 
\begin{eqnarray*}
\sum_{n=0}^{\infty}B_{p,n}(z) \frac{t^n}{[n]_p!}&=&\sum_{n=0}^{\infty}\left(\sum_{k=0}^{n}\pbinom{n}{k}B_{n-k} z^k\right)\frac{t^n}{[n]_p!}\\
&=&\sum_{k=0}^{\infty} \frac{(tz)^k}{[k]_p!}\sum_{n=k}^{\infty} \frac{B_{n-k}}{[n-k]_p!} t^{n-k}\\
&=&E(tz;p) \sum_{k=0}^{\infty}\frac{B_k}{[k]_p!} t^k.  
\end{eqnarray*}
So we  obtain 
\be\label{Eq_2}\sum_{n=0}^{\infty}B_{p,n}(z) \frac{t^n}{[n]_p!}=E(tz;p) \sum_{k=0}^{\infty}\frac{B_k}{[k]_p!} t^k.\ee
Now acting on the two sides of \eqref{Eq_2} by $T_p^1$ where we consider the two sides as  functions of $z$, i.e.  $t$ is fixed.
This gives 
\be\label{Eq_3}\sum_{n=0}^{\infty}B_{p,n}(z\dotplus 1) \frac{t^n}{[n]_q!}=E(t(z\dotplus 1);p) \sum_{k=0}^{\infty}\frac{B_k}{[k]_p!} t^k.\ee
Then substitute with $z=0$ and use  $B_n(0\dotplus 1)-B_n= 0$ if $n\neq 1$ and $1$  if $n=1$. This yields 
\be 
t+\sum_{k=0}^{\infty}\frac{B_k}{[k]_p!}=E(t(0\dotplus 1);p) \sum_{k=0}^{\infty}\frac{B_k}{[k]_p!}.
\ee
I.e. 
\be
\sum_{k=0}^{\infty}\frac{B_k}{[k]_p!}=\frac{t}{E(t(0\dotplus 1);p)-1},
\ee
and from \eqref{Eq_2}
\be \sum_{n=0}^{\infty}B_{p,n}(z) \frac{t^n}{[n]_p!}=\frac{tE_p(tz)}{E(t(0\dotplus 1);p)-1}.\ee

\end{proof}
\begin{exa}

Assume that $\delta_k=1$ for all $k\in\mathbb{N}_0$.  We have two cases:
\begin{enumerate} 
\item The case  $p=q$, $0<q<1$ gives 
\[E(tz;p)=e_q(tz),\quad E(t(z\dotplus 1);p)|_{z=0}=e_q(t),\] and we get 
\[\sum_{n=0}^{\infty}B_{p,n}(z)\frac{t^n}{[n]_p!}=\frac{te_q(tz)}{e_q(t)-1}.\]This yields that  $\left(B_{p,n}\right)_{n}$ is the set of $q$-Bernoulli polynomials introduced by Al-Salam in \cite[Eq. (5.1)]{Al-Salam59}.

\vskip .5 cm 

\item    The case $p>1$. In this case,  there exists $0<q<1$ such that $p =1/q$. Hence,  
\[E(tz;p)=E_{q}(tz),\quad E(t(z\dotplus 1);p)|_{z=0}=E_{q}(t),\] and we get 
\[\sum_{n=0}^{\infty}B_{p,n}(z)\frac{t^n}{[n]_p!}=\frac{tE_{q}(tz)}{E_{q}(t)-1}.\]This yields that
\[B_{p,n}(x)=q^{n\choose 2}b_n(x) \quad (n\in\mathbb{N}_0),\]
where $\left(b_n(x)\right)_n$ is the sequence of Bernoulli polynomials introduced by Al-Salam in \cite[Eq. (5.2)]{Al-Salam59}.
\end{enumerate}
\end{exa}

\begin{exa}

Let $\delta_k:=\frac{(-1;p)_k}{2^k}$ $(k\in\mathbb{N}_0)$.  If $p=q$, $0<q<1$, then 
\[E(tz;q)=e_q(tz),\quad E(t(z\dotplus 1);q)|_{z=0}=e_{q}(t/2)E_q(t/2),\]
and   
\[\sum_{n=0}^{\infty}B_{p,n}(z)\frac{t^n}{[n]_p!}=\frac{te_{
q}(tz)}{e_{q}(t/2)E_q(t/2)-1}.\]This yields that
\[B_{p,n}(x)=b_n(x;q)\quad (n\in\mathbb{N}_0),\]
where$ (b_n(x;q))_n$ is the sequence of $q$-Bernoulli polynomials  defined in \cite[Eq. (1.4)]{Ismail-Mansour}.  If $p>1$, then there exists $q$, $0<q<1$, such that $p=1/q$. Consequently,  
\[E(tz;p)=E_{q}(tz),\quad E(t(z\dotplus 1);p)|_{z=0}=e_{q}(t/2)E_{q}(t/2),\] and 
\[B_{p,n}(x)=q^{n \choose 2}\,B_n(x;q) \quad (n\in\mathbb{N}_0),\]
where$ (B_n(x;q))_n$ is the sequence of $q$-Bernoulli polynomials  defined in \cite[Eq. (2.1)]{Ismail-Mansour}.

\end{exa}
Given an entire function $f(z)=\sum_{n=0}^{\infty}a_n z^n$, the function $g(z)=\sum_{n=0}^{\infty}a_n\dfrac{B_{p,n+1}(z)}{[n+1]_p}$ satisfies the functional equation \eqref{FEP} in the domain of analyticity of the function $g(z)$, which depends on the large $n$ asymptotic of the polynomials $(B_{p,n}(x))_n$.  

In \cite{Mansour-Sahar}, the authors used the Darboux method to prove that  the $q$-Bernoulli polynomials $(B_n(x;q))_n$ ($0<q<1$),  we  introduced in \cite{Ismail-Mansour},  have the  larger $n$ asymptotics  
\be
\begin{gathered}
\dfrac{B_n(x;q)}{[n]_q!}=-2\dfrac{\cos \frac{n\pi}{2}}{(2\xi_1)^{n+1}}\dfrac{\text{Cos}_q 2\xi_1 x\, \text{Cos}_q \xi_1}{\text{Sin}'_q(\xi_1)}\\
-2\dfrac{\sin \frac{n\pi}{2}}{(2\xi_1)^{n+1}}\dfrac{\text{Sin}_q 2\xi_1 x\, \text{Cos}_q \xi_1}{\text{Sin}'_q(\xi_1)}+o(r^n),
\end{gathered}
\ee
where $r< 2\xi_1$,  and $\xi_1$ is the smallest positive zero of $\text{Sin}_q x$.  Recall that the functions $\text{Sin}_q x$ and $\text{Cos}_qx$ are defined  by the  $q$-Euler's formula
\[E_q(ix)=\text{Cos}_qx +i \text{Sin}_q x. \]
It is known that the zeros of $\text{Sin}_qx$ and $\text{Cos}_qx$ are real simple interlacing zeros. See. \cite{Ism1}. 
Hence, there exists $M>0 $ such that 
\be\label{Ineq:1}
\left|\dfrac{B_n(x;q)}{[n]_q!}\right|\leq \frac{M}{(2\xi_1)^n }E_q(2\xi_1 |x|) \quad (x\in\mathbb{C},\, n\in\mathbb{N}),
\ee
and 
\be\label{Ineq:2}
\left|\dfrac{B_n(q)}{[n]_q!}\right|\leq \frac{M}{(2\xi_1)^n } \quad ( n\in\mathbb{N}).
\ee	
This leads to the following theorem. 
\begin{thm}
Let $p>1$ and  $f(z):=\sum_{n=0}^{\infty}a_n z^n$ be an entire function of $p$-exponential growth of order $c\ln p $, $c<1$.  There exists an entire function $g(z)$ such that 
\be \label{FE} g(z\dotplus 1)-g(z)=f(z), \ee
where $g(z\dotplus 1)=T_p^1 g(z)$ , and the operator $T_p^1$ is the one defined in \eqref{T:def} with $\delta_k:=\frac{(-1;p)_k}{2^k}$.

\end{thm}
 
\begin{proof}Since $p>1$, then there exists  $q$, $0<q<1$ such that   $p=1/q$.   
Set  \[g(z):=\sum_{n=0}^{\infty}\frac{a_n}{[n+1]_q} q^{-\frac{n(n-1)}{2}}B_{n+1}(z;q) .\] Since the polynomials, $(B_n(z;q))_n$, satisfy  the identity 
\[B_n(z\dotplus 1;q)-B_n(z;q)=q^{\frac{(n-1)(n-2)}{2}}[n]_qz^{n-1}\quad (n\in\mathbb{N}), \]
then one can verify that the function $g$ satisfies \eqref{FE}. Since $f$ is of $q^{-1}$-exponential growth of order $-c\ln q$,    $c<1$,  and a certain type $\alpha$.  Then there exists $K>0$
 such that $|a_n|\leq K q^{\frac{(n-\alpha)^2}{2c}}$, $n\in\mathbb{N}$. Consequently, from \eqref{Ineq:1},
 \[\left|\sum_{n=0}^{\infty}\frac{a_n}{[n+1]_q} q^{-\frac{n(n-1)}{2}}B_{n+1}(z;q)\right|\leq KM E_q(2\xi_1 |z|)\sum_{n=0}^{\infty}q^{\frac{(n-\alpha)^2}{2c}}q^{-\frac{n(n-1)}{2}}\frac{[n]_q!}{(2\xi_1)^{n+1}} <\infty. \]
 Hence $g(z)$ is an entire function. 

\end{proof}

\end{document}